\font\Large=cmr10 scaled\magstep2  
\hfuzz=4pt

\magnification=\magstep1

\input amssym.def 
\input amssym.tex

\def\A{{\Bbb A}} \def\C{{\Bbb C}} \def\R{{\Bbb R}} \def\Z{{\Bbb Z}} \def\P{{\Bbb 
P}} \def\ra{{\longrightarrow}}

\def\O{{\cal O}} \def\ot{{\otimes}} \def\lb{{\langle}} \def\rb{{\rangle}} 
\def\E{{\cal E}} \def\F{{\cal F}}
\def\U{{\cal U}} \def\llongrightarrow{ {-\!\!-\!\!\!\longrightarrow}} 
\def\longdownarrow{\Big\downarrow }

\centerline{\Large Noncommutative smooth spaces}\par

\bigskip

\medskip

\centerline {{\bf Maxim Kontsevich, Alexander Rosenberg}} \vskip 1.5truecm

\bigskip

\centerline{\bf 0.  \rm CONVENTIONS AND NOTATIONS}\medskip

We will work in the category $Alg_k$ of associative unital algebras over a fixed 
base field $k$.  If $A\in
Ob(Alg_k)$, we denote by $1_A\in A$ the unit in $A$ and by $m_A:A\otimes A\ra A$ 
the product.  For an algebra $A$, we
denote by $A^{opp}$ the opposite algebra, i.e.  the same vector space as $A$ 
endowed with the multiplication
$m_{A^{opp}}(a\otimes b):=m_A(b\otimes a)$.  If $A$ and $B$ are two algebras, 
then 
$A\otimes _kB$ is again an
algebra.  Also, $A\star B$ denotes the free product of $A$ and $B$ over $k$, the 
coproduct in the category $Alg_k$.  By
$A$-$mod$ we denote the abelian category of left $A$-modules.  Analogously, 
$mod$-$A$ are right modules (the same
as $A^{opp}$-modules) and $A$-$mod$-$A$ are bimodules over $k$, or, 
equivalently, 
$A\otimes _kA^{opp}$-modules. We shall write $\otimes $ instead of $\otimes _k$. 
For a vector space $V$, we denote by $Sym^*(V)$ and $\otimes ^*(V)$ resp.  the 
free commutative associative
(polynomial) and free associative (tensor) $k$-algebra respectively, generated 
by 
$V$.

\bigskip

\centerline {\bf 1.  \rm SMOOTH ALGEBRAS}

\medskip

A typical example of an algebra for this text is a free finitely generated 
algebra 
$k\lb x_1,\dots, x_d\rb$, in
 contrast to usual noncommutative algebras close to commutative ones, like 
universal enveloping algebras,
 algebras of algebraic differential operators etc.
  We consider the free algebra  as ``the algebra of functions
  on the non-commutative affine space'' which we denote $N\A^d$. Recently there 
were several
  attempts
  to understand  algebraic geometry  of this space.  Gelfand and Retakh started 
the study
 of basic identities in the free skew field with $d$ generators (see [GR]), 
which is  can be considered
  as the  ``generic point'' in $N\A^d$.  One of us (see [K])
   tried to develop ``differential geometry''
    in the algebra of noncommutative formal power series, i.e. at the formal 
neighborhood
     of $N\A^d$ at zero.
     Kapranov (see [Ka]) 
   described a differential-geometric
   picture related with the completion of the free algebra with respect to the 
commutator filtration,
    i.e. with the
    algebra of functions on the infinitesimal neighborhood of the usual affine 
space $\A^d$  in  $N\A^d$. 
    Here we would like to study $N\A^d$ as the whole space, without completions 
and localizations. 
    The basic intrinsic property of the free algebra is {\it smoothness. \rm} 

\medskip

\noindent {\bf 1.1.  Equivalent definitions}

\medskip

Few years ago J.\ Cuntz and D.\ Quillen gave the following

\medskip 

\noindent {\bf Definition}  ([CQ1], Definition 3.3 and Proposition 6.1).  {\it 
An 
algebra $A$ is quasi-free (= formally smooth) iff
it satisfies one of the following equivalent properties:

 1) (Lifting property for nilpotent extensions) for any algebra $B$, a two-sided 
nilpotent ideal $I\subset B$
  ($I=BIB,\,\,\,I^n=0$ for $n\gg 0$), and for any algebra homomorphism $f:A\ra 
B/I$, there exists an algebra
  homomorphism ${\tilde f}:A\ra B$ such that $f=pr_{B\ra B/I}\circ{\tilde f}$, 
where $pr_{B\ra B/I}:B\ra B/I$ is
  the natural projection.

2) $Ext^2_{A-mod-A} (A,M)=0$ for any bimodule $M\in Ob (A$-$mod$-$A)$.

3) The $A$-bimodule $\Omega ^1_A:=Ker(m_A:A\otimes A\ra A)$ (see 1.1.2) is 
projective.\rm \par

\medskip

The definition of formal smoothness via the lifting property 1) is analogous to 
the Grothendieck's (actually Quillen's) definition
 of formally smooth algebras in the commutative case.  The equivalence of 
properties 1)-2)-3) is an easy exercise in
 homological algebra.

Properties 1),2) are not constructive: one can not check them directly.  The 
property 3) looks better, although it
is still not clear {\it a priori } how to work with it in complicated 
infinite-dimensional examples.  J.~Cuntz and
D.~Quillen have found another characterization of smooth algebras which is 
convenient for calculations in
practice.  In order to give this characterization we need an auxiliary 
definition.

Let $TA$ be the algebra generated by symbols $\underline{a}, \underline{Da}$, 
where $a\in A$, subject to the
following relations:

a) the map $a\mapsto \underline{a} $ is a homomorphism of unital $k$-algebras,

b) the map $a\mapsto \underline {Da}$ is $k$-linear,

      c) $\underline{D(a\cdot b)}=\underline{a}\cdot 
\underline{Db}+\underline{Da}\cdot
\underline{b}$, (the Leibnitz rule).

Algebra $TA$ is naturally $\Z_{\ge 0}$-graded:  $deg(\underline{a})=0,\,\,\, 
deg(\underline{Da})=+1$.  We identify
$A$ with the subalgebra of $TA$ consisting of elements of degree zero.

\medskip \noindent {\bf 1.1.2.  Differential envelope} \medskip Before going 
further on we would like to notice
that $TA$ is isomorphic as an abstract $\Z$-graded algebra to the so-called {\it 
differential envelope $\Omega
A$}, the universal differential $\Z$-graded super-algebra containing $A$.

By definition, $\Omega A$ is the universal $\Z$-graded super-algebra containing 
$A$ and endowed with an odd
differential $d$ of degree $+1$ (satisfying the Leibnitz rule in super sense 
$d(a\cdot b)=d(a)\cdot b+(-1)^{deg
(a)} a\cdot d(b)$) and such that $d^2=0$.  One can see immediately that $d1=0$.  
As a vector
 space, the $n$-th graded
component $(\Omega A)^n$ for $n\ge 0$ is isomorphic to $A\otimes (A/k\cdot 
1_A)^{\otimes n}$.  The isomorphism is
given by the following map $$a_0\otimes a_1\otimes\dots\otimes a_n\longmapsto 
a_0\cdot da_1\cdot\dots\cdot da_n\in (\Omega
A)^n\,\,\,.$$ The space of $1$-differentials $\Omega^1_A:=(\Omega A)^1$ can be 
furthermore identified via the map
$$A^{\otimes 3}\ra \Omega^1_A,\,\,\,a_1\otimes a_2\otimes a_3 \longmapsto 
a_a\cdot 
da_2\cdot a_3$$ with the 
quotient space $A^{\otimes 3}/\partial_4 (A^{\otimes 4})$, where $$\partial_4 
(a_1\otimes a_2\otimes a_3\otimes
a_4)=a_1 a_2\otimes a_3\otimes a_4 -a_1\otimes a_2 a_3\otimes a_4+ a_1\otimes 
a_2\otimes a_3 a_4)\,\,\,.$$

The map $\partial_4$ can be, evidently, included into one of standard complexes 
in 
homological algebra
 $$\dots\ra A^{\ot 4}\ra A^{\ot 3}\ra A^{\ot 2} \ra A\ra 0$$ This complex can be 
contracted using the unit element
 $1_A\in A$, and thus has vanishing cohomology.  Therefore 
$Cokernel(\partial_4)= 
Kernel(\partial_2)$ and we get
 an equivalent description of $\Omega^1_A$ as the kernel of the multiplication 
map 
$A\otimes A\ra A$.

If one considers analogous definitions in the purely commutative case, one can 
see 
immediately the difference
 between $TA$ and $\Omega A$.  In the case of $A$ equal to the algebra of 
functions on a smooth affine variety
 $X$, the commutative analogue of the algebra $TA$ is the algebra of functions 
on 
the total space $TX$ of the
 tangent bundle to $X$.  On the contrary, the commutative analogue of $\Omega A 
$ 
is the algebra of differential
 forms on $X$, which is the same as the algebra of functions on supermanifold 
$\Pi 
TX$, the total space of the odd
 tangent bundle to $Z$ (as usual in super-mathematics,  letter $\Pi$ denotes the 
functor of changing the
 parity).  

\medskip 

\noindent {\bf 1.1.3.  Another criterion of formal smoothness} 

\medskip

\proclaim {Theorem}.  An algebra $A$ is formally smooth iff it satisfies the 
property

\it 4) there exists a derivation $\underline{D}:TA\ra TA$ of degree $+1$ such 
that
  $\underline{D}(\underline{a})=\underline{Da}$ for all $a\in A$.  \rm \par

If $A$ is an algebra generated by elements $(a_i)_{i\in I}$, then in order to 
define $\underline{D}$ one should 
 define only elements $\underline{D}(\underline{Da_i})\in TA$ of degree $+2$ 
satisfying certain relations.  In the
 case when the generating set $I$ is finite and all relations between $(a_i)$ 
follow from a finite number of
 relations we say that $A$ is finitely presented.  It is easy to see that in 
this 
case algebra $TA$ is also
 finitely presented and in order to check that  some candidates for 
$\underline{D}(\underline{Da_i})$ satisfy
 needed relations one should make only a finite calculation.

 The theorem analogous to the one above holds in the commutative case (if the 
characteristic of the ground field
$k$ is zero).  The derivation $\underline{D}$ in commutative case is a vector 
field on the total space of the
tangent bundle.  One can easily see that choices of $\underline{D}$ correspond 
to 
symmetric connections (i.e.
connections with vanishing torsion) on the tangent bundle to $Spec(A)$.

\medskip

\noindent {\bf 1.1.4.  Homological properties of modules over smooth algebras}

\medskip

\proclaim {Theorem}. {\rm ([CQ1], Proposition 5.1)} If $A$ is a formally smooth 
algebra, then $A$ is hereditary, i.e.  the category $A$-mod has
	homological dimension $\le 1$.  In other words, any submodule of a 
projective module is projective.
	Equivalently, every $A$-module admits a projective resolution $P_{-1}\ra 
P_0$ of length $2$.  \par

The proof is immediate, because $$Ext^n_{A-mod}(\E,\F)= Ext_{A-mod-A}^n(A, 
Hom_{k-mod}(\E,\F))\,\,\,\,\forall
\E,\F\in Ob(A-mod)\,\,\,.$$

\proclaim {Definition}. We call an algebra $A$ smooth, if it is formally smooth 
and finitely generated.

\proclaim {Theorem}.  If $A$ is a smooth algebra and $\E,\F$ are  $A$-modules 
 finite-dimensional over $k$, then the vector spaces 
$Ext^0_{A-mod}(\E,\F)=Hom_{A-mod}(\E,\F)$ and
 $Ext^1_{A-mod}(\E,\F)$ are both finite-dimensional.

\par

The statement of the theorem is evident for $Ext^0$.  The space $Ext^1$ 
coincides 
with the set of equivalence
 classes of structures of an $A$-module on $\E\oplus \F$ such that natural maps 
$$\F\ra \E\oplus \F\ra \E$$ are
 morphisms of $A$-modules, modulo the action of the vector space 
$Hom_{k-mod}(\E,\F)$ considered as an abelian
 group.  If $A$ is finitely generated, then the set of $A$-module structures on 
$\E\oplus \F$ as above is a
 finite-dimensional vector space over $k$, and we get the statement of the 
theorem.

\medskip \noindent {\bf 1.2.  Examples of smooth algebras} \medskip

E1) The free algebra $k\lb x_1,\dots,x_d\rb=\otimes^*(k^d)\,\,\,d\ge 0$.  We 
think 
of this algebra as the one 
corresponding to a
noncommutative affine space $N\A^d_k$ (see Sect.  3 where we define 
noncommutative 
spaces in general).

E2) The matrix algebra $Mat(n\times n,k)$.

E3) The algebra of upper-triangular matrices $$A=\{(a_{ij})_{1\le i,j\le n} 
|\,\,a_{ij}\in
   k,\,\,\,a_{ij}=0\,\,\,\,\,for\,\,\,\,\, i>j\}\subset Mat(n\times 
n,k)\,\,\,,$$

E4) $\O(C)$, the algebra of functions on a smooth affine curve $C$ over $k$.

E5) $Paths(\Gamma)$, the algebra of finite paths in a finite oriented graph 
$\Gamma$.  The basis of 
 $Paths(\Gamma)$ consists of sequences 
$(v_0,e_{0,1},v_1,e_{1,2},\dots,e_{n-1,n}, 
v_n),\,\,\,n\ge 0$ of vertices
 $v_i$ of $\Gamma$ and oriented edges $e_{i,i+1}$ connecting $v_i$ with 
$v_{i+1}$ 
(there could be several edges
 connecting two given vertices).  The unit element in $Paths(\Gamma)$ is equal 
to 
the sum over vertices $v$ of
 $\Gamma$ of paths of zero length $(v)$.  The multiplication in $Paths(\Gamma)$ 
is 
given by concatenation of paths:
 $$(v_0,e_{0,1},\dots,v_n)\cdot (v'_0,e'_{0,1},\dots,v'_{n'})=$$
 $$(v_0,e_{0,1},\dots,v_n,e'_{0,1},\dots,v'_{n'})$$ if $v_n=v'_0$, and $0$ 
otherwise.

The last example E5) contains E1) and E3) as particular cases.  In E1) the 
corresponding graph has one vertex and
$d$ loops, and in E3) it is the chain of $n-1$ consecutive oriented 
edges. 

 There are several constructions using which one can
construct new smooth algebras from the old ones.  In the following list 
$A,A_1,A_2$ denote smooth algebras, and
the result of the constructions is automatically a smooth algebra:

C1) $A_1\oplus A_2$, the direct sum of two smooth algebras (or an arbitrary 
finite 
number of smooth algebras),

C2) $A_1\star A_2$, the free product of smooth algebras.

C3) (Localization) If $S$ is a subset of a formally smooth algebra $A$, then the 
algebra $A[S^{-1}]$ obtained by formally adjoining the inverses of elements of 
$S$ 
is formally smooth. In particular, if $A$ is smooth and $S\subset A$ is an 
arbitrary finitely generated multiplicative subset of $A$, then $S^{-1}A$ is 
smooth.  More generally, one can invert not only individual elements of $A$ but 
square matrices with coefficients in $A$ (i.e.
$S$ could be a subset of $\coprod _{n\ge 1} Mat(n\times n,A)$), or even 
rectangular matrices.

C4) $End_{A-mod} (P)$, where $P$ is a finitely-generated projective $A$-module.

C5) (The total space of the tangent bundle) $TA$, the algebra defined in 1.1.

\medskip 

\noindent {\bf 1.3.  Representation spaces of smooth algebras} 

\medskip

Let $A$ be a smooth $k$-algebra.  We associate with $A$ an infinite sequence
$(Repr^A_n,\,\,\,n\ge 1)$ of smooth affine varieties over $k$.  The $n$-th 
variety 
$Repr^A_n$ parameterizes
homomorphisms from $A$ to the standard matrix algebra $Mat(n\times n, k)$.  The 
set of $B$-points of it, where $B$
is a {\it commutative \rm} algebra over $k$, is defined as 
$$Repr^A_n(B):=Hom_{Alg_k}(A,Mat(n\times
n,B))\,\,\,.$$ 

It is clear that since $A$ is finitely-generated, the scheme $Repr^A_n$ is an 
affine scheme of finite
type over $k$.  Moreover, from the definitions of smoothness in both commutative 
and noncommutative case, it
follows that the scheme $Repr^A_n$ is smooth.  The affine algebraic group 
$PGL_k(n)$ acts via conjugations on
$Repr^A_n$.

Every element $a\in A$ gives tautologically a matrix-valued function 
$\hat{a}\in\O(Repr^A_n)\otimes Mat(n\times
n,k)$ on $Repr^A_n$.

\medskip

\noindent {\bf 1.3.1.  Examples}

\medskip

For $A=k\lb x_1,\dots, x_d\rb$, the representation scheme $Repr^A_n$ coincides 
with the affine space $\A^{d n^2}$
over $k$.

For $A=k\oplus k=k\lb p\rb/(p^2=p)$, the scheme $Repr^A_n$ is the disjoint union 
over all integers $m,\,\,\,0\le
 m\le n$ of certain bundles over Grassmanians $Gr(m,n)$.  Fibers of these 
bundles 
are affine spaces parallel to
 the fibers of the cotangent bundles.  The component of $Repr^A_n$ corresponding 
to $m$ (--- the rank of the image
 of $p\in A$) has dimension $2m(n-m)$.

The example of $A=k\oplus k$ has an interesting feature.  One can see that 
$Repr^A_n$ carries a natural
non-degenerate closed $2$-form $$\omega=Trace(\hat{p}\, d\hat{p}\wedge  
d\hat{p})$$ 
where $\hat{p}$ is the
matrix-valued function on $Repr^A_n$ corresponding to the generator $p$ of $A$.  

\medskip 

\noindent {\bf 1.3.2. Noncommutative analogues of differential-geometric  
structures} 

\medskip

We would like to propose the following general principle which could help to 
look 
for noncommutative versions of
usual geometric notions (see also [K]).

\it A noncommutative structure of some kind on $A$ should give an analogous 
``commutative'' structure on all schemes
 $Repr^A_n,\,\,\,n\ge1$.  \rm

Here we make a list of several natural candidates:

{\bf Functions}.  Elements of the vector space $A/[A,A]$ give rise to {\it 
functions \rm} on $Repr^A_n$.  Namely,
 any $a\in A$ modulo commutators $[A,A]$ gives the function $Trace(\hat{a})$.  
Thus, we get a homomorphism of
 vector spaces $$A/[A,A]\ra \O(Repr^A_n)\,\,\,.$$ It has an obvious extension to 
the homomorphism of algebras
 $$Sym^*(A/[A,A])\ra \O(Repr^A_n)\,\,\,.$$

We will denote the element of $\O(Repr^A_n)$ corresponding to $\phi\in 
Sym^*(A/[A,A])$ by $Trace(\phi)$.

In the following examples we can multiply our differential-geometric objects by 
an 
arbitrary element of the
 algebra $Sym^*(A/[A,A])$.

{\bf Vector fields}.  Any derivation $\xi\in Der(A)$ of $A$ gives a vector field 
$\tilde {\xi}$ on $Repr^A_n$.

{\bf Differential forms}.  Any element of $\Omega A$ (see 1.1.2) gives a  
matrix-valued differential form on
 $Repr^A_n$, and thus a scalar-valued differential form after taking the trace 
in 
the matrix algebra.  Again,
 using the vanishing of the trace of commutators one can see that the map 
$\Omega 
A \ra \Omega^*(Repr^A_n)$ goes
 through the quotient space $$\Omega A/[\Omega A,\Omega A]_{super}\ra 
\Omega^*(Repr^A_n)$$ where the commutator in
 superalgebra $\Omega A$ is understood in super sense.  Also, one can easily see 
that the differential $d$ in
 $\Omega_A$ maps to the de Rham differential in usual forms 
$\Omega^*(Repr^A_n)$.  
For example, algebra $A=k\oplus
 k$ has a closed $2$-form.  The definition of ``noncommutative differential 
forms'' as elements of $\Omega
 A/[\Omega A,\Omega A]_{super}$ is due to M.~Karoubi.

{\bf De Rham cohomology}.  It is natural in the light of the previous discussion 
to define the de Rham cohomology of
 $A$ as the cohomology of the complex $\Omega A/[\Omega A,\Omega A]$.  Using 
results of Cuntz and Quillen [CQ2], one can
 show that $H^n_{dR}(A)$ defined in this way coincides with the $Z/2\Z$-graded 
periodic cyclic homology 
 $HP_n(A)$ for $n=0,1$, and with the reduced periodic cyclic homology
 
$$\widetilde{HP}_{n(mod\,\,\,2)}(A):={HP}_{n(mod\,\,\,2)}(A)/{HP}_{n(mod\,\,\,2)
}(
k)$$ for $n\ge 2$.  In other
 words, $H^n_{dR}(A)$ coincides with $HP_{n(mod\,\,\,2)}(A)$ for 
$n=0,1,3,5,\dots$ 
and has one dimension less than
 $HP_{n(mod\,\,\,2)}(A)$ for $n=2,4,6,\dots$.

{\bf Volume element, polyvector fields}.  It looks plausible from examples that 
one can define noncommutative
 analogues of polyvector fields and volume elements, but we would like to 
postpone 
the discussion of these notions
 to the future.  It seems that  the divergence of a derivation of $A$ with 
respect 
to a volume element belongs
  to the
 symmetric square of the vector space $A/[A,A]$.  As an exercise to the reader 
we 
leave the following simple

\proclaim {Lemma}.  If $A=k\lb x_1,\dots x_d\rb$ then there exists a linear map 
$$div:  Der(A)\ra Sym^2(A/[A,A])$$
 such that for any $\xi\in Der(A)$ the divergence of the corresponding vector 
field $\tilde{\xi}$ on
 $\A^{dn^2}=Repr_n^A$ with respect to the standard volume element, is equal to 
$Trace(div(\xi))$.  \par

   The last remark which we would like to make is that the correspondence ({\it 
Smooth algebras\rm}) $\ra$ ({\it
Representation schemes\rm}) gives a ``justification'' of the notation $TA$ in 
1.1. 
 Namely, there is a natural
isomorphism between $Repr_n^{TA}$ and the total space of the tangent bundle 
$TRepr_n^A$.  

\medskip 

\noindent {\bf 1.3.3.  Possible relations to matrix integrals and to M-theory} 

\medskip

 One can imagine that in the case $k=\C$, for a given ``volume element'' $vol$ 
on 
$A$, a ``function'' $f\in
Sym^*(A/[A,A])$ and for a the set of real points $\gamma_n$ in $Hom(A, 
Mat(n\times 
n,\C))$ under some
anti-holomorphic involution, one can take the integral of $exp(Trace(f))\times 
vol$ and get a ``matrix model'', an
infinite sequence of numbers $$I_n:=\int_{\gamma_n} e^{Trace(f)} vol$$  
parameterized by the dimension
$n=1,2,\dots$.  In mathematical physics such integrals were extensively studied.  
Typically, one integrates over
the space of hermitean $n\times n$ matrices, or over real, or unitary matrices.  
In multi-matrix models the
integration is taken over the set of, say, $d$-tuples $(X_1,\dots, X_d)$ of 
hermitean matrices.  It is believed
that the asymptotic behavior of $I_n$ as $n\ra \infty$ is related with some kind 
of string theory.

Also, one of recently proposed matrix models, so called $M$-theory, is 
formulated in the same fashion.  Roughly
speaking, $M$-theory on the space-time manifold $X=\R^d$ is a matrix theory 
corresponding to the free algebra with
$d$ generators, the noncommutative affine space.  $M$-theory on  curved spaces 
 should correspond to nontrivial smooth noncommutative algebras.

\medskip 

\noindent {\bf 1.3.4.  Double tangent space and formal noncommutative structure} 

\medskip

Here we would like to give some examples of natural non-classical structures on 
manifolds $Repr_n^A$.  First of
  all, there is a natural vector bundle $T_{(2)}$ on the square 
$Repr_n^A\times_k 
Repr_n^A$ together with the
  identification of its pullback via the diagonal embedding $\Delta :  
Repr_n^A\ra 
Repr_n^A \times_k Repr_n^A$
  with the tangent bundle $T_{Repr_n^A}$.

 In what follows we will describe the bundle $T_{(2)}$.  First of all, 
$k$-points 
of $Repr_n(A)$ can be identified
with equivalence classes of pairs $(\E,pr_\E)$ where $\E$ is an $A$-module 
finite-dimensional over $k$, and
$pr_\E:A^n\ra \E$ is an epimorphism to $\E$ from the standard $n$-dimensional 
free 
$A$-module such that the set of
$n$ canonical generators of $A^n$ maps to a basis of $\E$ over $k$.  The fiber 
of 
$T_{(2)}$ at the pair
$\bigl((\E,pr_\E),(\F,pr_\F)\bigr)$ is defined as 
$$Hom_{A-mod}(Ker(pr_\E),\F)\,\,\,.$$ Let us prove that it is a
finite-dimensional vector space.  It is easy to see that it fits into an exact 
sequence $$0\ra
Hom_{A-mod}(\E,\F)\ra\F^n\ra Hom_{A-mod}(Ker(pr_\E),\F)\ra Ext^1_{A-mod}(E,F)\ra 
0\,\,\,$$ All spaces 
except the one in question are
finite-dimensional by the theorem 1.1.4.  Also, it is easy to see that on the 
diagonal $\E=\F,\,\,pr_\E=pr_\F$, the
space $Hom_{A-mod}(Ker(pr_\E),\F)$ coincides with the tangent space to 
$Repr_n^A$ 
at the point $(\E,pr_\E)$.

Another structure on $Repr_n^A$ is the {\it formal noncommutative structure\rm} 
in 
the sense of Kapranov.  There
 is a canonical sheaf $\O^{noncomm}_{Repr^A_n}$ of formal noncommutative  
functions 
on $Repr_n^A$ with the
 quotient $\O_{Repr_n^A}$.  For example, if $A=k\lb x_1,\dots, x_d\rb$ and $n=1$ 
then the global sections of the
 sheaf $\O^{noncomm}_{Repr^A_n}$ is the projective limit ${lim} (A/I_n)$ where 
$I_n, n\ge 1$ is a decreasing
 sequence of two-sided ideals $$I_n=\sum_{l\ge 1,(m_1,\dots,m_l):\sum m_i=n} 
A\cdot A^{[m_1]}\cdot A\cdot
 A^{[m_2]}\cdot\dots\cdot A^{[m_l]}\cdot A\,\,\,.$$ Here $A^{[m]}$ for $m\ge 1$ 
is 
defined as the linear span of
 the set of commutators of depth $m$:  $$[a_0,[a_1,[\dots[a_{m-1},a_m]\dots]\in 
A,\,\,\,a_0,a_1,\dots,a_m\in
 A\,\,\,.$$ We refer the reader to [Ka] for the definition of the formal  
noncommutative structure and its
 differential-geometric meaning.  We would like only to mention that Kapranov 
proposed a general construction of
 formal noncommutative structure on moduli spaces of objects in abelian 
categories, like the category of coherent
 sheaves on algebraic varieties etc.  His construction admits a useful 
extension.  
Namely, if $\cal C$ is a
 $k$-linear triangulated category and $\E_0$ is a fixed object, then one can try 
to consider the ``moduli space''
 of pairs $(\E,p)$, where $\E$ is an object of $\cal C$ and $p:\E_0\ra \E$ is a 
morphism.  The tangent complex for
   such a pair is $RHom_{\cal C} (\E, Cone(p))$.  It carries a natural structure 
of an associative non-unital
 algebra, and the arguments from [Ka] are applicable there.

\medskip

\noindent {\bf 1.4.  What we would like to do?}

\medskip

It is clear from all the previous discussion that the subject of smooth 
noncommutative geometry merits further
development.  All our previous constructions and examples are {\it affine\rm}, 
in 
particular they give affine
representation schemes.  It would be desirable to give a definition of ``smooth 
noncommutative schemes'', and
also of ``quasi-projective noncommutative schemes'', such that it contains 
finitely generated (or maybe finitely
presented) smooth algebras as an affine case, and give rise to smooth schemes 
via 
certain functor of
``representations to $Mat(n\times n, k)$''.  We expect that representation 
schemes 
carry bitangent bundles,
Kapranov's formal noncommutative structure, etc.

  In the second part of this text we describe a general approach to 
noncommutative 
algebraic geometry based on flat
 topology.  We choose (temporarily) 
 the name ``spaces'' (or ``noncommutative spaces'') instead of ``schemes'' for 
several reasons.
  Noncommutative schemes were introduced in earlier works of one of us (see 
[R2]), 
where 
analogues
 of Zariski topology were studied. 
 Also, our formalism describes not only usual schemes but  algebraic spaces too.
   As the reader will see,
  an amazing variety 
  of basic constructions in usual algebraic geometry can be extended to the 
noncommutative setting.
  Our category of noncommutative spaces could be used in other situations. For 
example,  quantum
 projective spaces arising from Sklyanin algebras also fit into our definition. 
 
  In the third part we study one
 particular space which we call ``noncommutative projective space'' and denote 
by 
$N{\P}^n_k$.  We beleive
 that $N{\P}^n_k$ is one of the principal examples of what should be called 
``smooth projective noncommutative
 variety''.

\bigskip

\centerline {\bf 2.  \rm NONCOMMUTATIVE SPACES}

\medskip

\noindent {\bf 2.1.  Covers and refinements}

\medskip

Here we describe  intermediate objects which are not yet  spaces (there are
 no structure sheaves on them). These objects are machines producing abelian 
categories. 
Essentially all definitions can be made in general monoidal categories (even in 
non-additive categories), instead of 
the category
 of vector spaces over $k$ with the monoidal structure given by the tensor 
product.
  
\medskip

\noindent {\bf 2.1.1.  Category of finite covers}

\medskip

\proclaim {Definition}.  Objects of category $Covers_k$ are given by the  
following

\it DATA:

1) an associative algebra $B\in Ob(Alg_k)$,

2) $M\in Ob(B-mod-B)$, a bimodule over $B$,

3) $m_M:M\ra M\otimes_B M$, a homomorphism of bimodules,

4) $e_M:M\ra B$, also a homomorphism of bimodules,

satisfying the following

AXIOMS:

1) $M$ is faithfully flat as the right $B$-module, i.e.  the functor 
$M\otimes_B:  
B-mod \ra B-mod$ is exact and
 does not kill any non-zero morphism,

2) $M$ with $m_M$ and $e_M$ is a coassociative coalgebra with counit in the 
monoidal category of $B\otimes
 B^{opp}$-modules.  \rm \par

\medskip

We will usually denote objects of the category $Covers_k$ as pairs $(B,M)$ 
skipping data
 $m_M$ and $e_M$.
 Morally, we consider pairs $(B,M)$ as ``noncommutative stacks'' together with a 
finite affine cover with affine
 pairwise intersections of members of the cover (see examples in the next 
subsection).

For every object $(B,M)$ of $Covers_k$ we have an abelian category $QCoh(B,M)$ 
whose objects are pairs $(\E,m_\E)$
 where $\E$ is a $B$-module and $m_\E:\E\ra M\otimes_B \E$ is a homomorphism of 
$B$-modules which defines a
 coaction of coalgebra $M$ on $\E$.  We call objects of $QCoh(B,M)$ \it 
quasi-coherent sheaves \rm on the
 ``noncommutative stack'' corresponding to $(B,M)$. One can justify the name 
``sheaf'' introducing
  an appropriate Grothendieck topology.

\proclaim {Definition}.  Morphisms $f$ in the category $Covers_k$ from 
$(B_1,M_1)$ 
to $(B_2,M_2)$ are given by the
 following

{\it DATA:

1) $f_B:B_2\ra B_1$, a morphism of algebras,

2) $f_M:M_2\ra M_1$, a morphism of $k$-vector spaces

satisfying the following

AXIOMS:

1) the diagram where vertical arrows are structure morphisms of $M_i$ as  
$B_i$-bimodules ($i=1,2$) is commutative:
 $$\matrix{\hfill B_2\otimes M_2\otimes B_2 & \buildrel f_B\otimes f_M\otimes 
f_B 
\over\llongrightarrow &
 B_1\otimes M_1\otimes B_1 \hfill \cr \longdownarrow \quad & & \quad 
\longdownarrow \cr M_2\quad & \buildrel f_M
 \over\llongrightarrow & \quad M_1 \cr}$$

2),3)
 two analogous diagrams including coproduct morphisms $m_{M_i}$ and counit 
morphisms $e_{M_i}$ respectively, are
 commutative.}

\medskip

Notice that the direction of the morphism $f=(f_B,f_M)$ is {\it opposite \rm} to 
the direction of the pullback
 morphisms of algebras and bimodules.

Every morphism $f=(f_B,f_M):(B_1,M_1)\ra (B_2,M_2)$ of covers defines a functor 
$f^*:QCoh(B_2,M_2)\ra
QCoh(B_1,M_1)$, the pullback of quasi-coherent sheaves.  This functor maps an 
object $(\E,m_\E)$ to the object
$(B_1\otimes_{B_2}\E,m'_\E)$ where the coaction morphism 
$$m'_\E:B_1\otimes_{B_2}\E\ra M_1\otimes_{B_1}
B_1\otimes_{B_2}\E=M_1\otimes_{B_2}\E$$ is defined in a natural way using $f_B$, 
$f_M$ and $m_\E$.

 \medskip 

 \noindent {\bf 2.1.2.  Examples of finite covers}\medskip

FC1) Let $S/Spec(k)$ be a separated quasi-compact scheme, e.g.  a 
quasi-projective 
scheme.  We choose a
finite cover $(\U_i)_{i\in I}$ of $S$ in Zariski topology by affine schemes.  It 
follows from separatedness that
pairwise intersections $\U_i\cap \U_j$ are again affine.  The algebra $B$ and 
the 
bimodule $M$ are defined as
$$B:=\O(\sqcup_i \U_i),\,\,\,\,M:=\O(\sqcup_{i,j} ( \U_i\cap \U_j))$$ and the 
structure of coalgebra on $M$ is the
natural one.  It follows from the usual descent theory that the category  
$QCoh(B,M)$ is equivalent to the category of
quasi-coherent sheaves on $S$.

FC1') The same statement holds for affine covers of separated quasi-compact 
schemes in fpqc topology, for
algebraic stacks etc.

FC2) Let $A$ be an associative algebra and $M:=A$ with the natural structure of 
an 
$A$-bimodule, and a coalgebra
 in the category $A$-$mod$-$A$.  Then $QCoh(A,M)$ is equivalent via the 
tautological functor to the category
 $A$-$mod$.

FC3) Let $A$ be an associative algebra, and $B\supset A$ be a larger algebra 
such 
that $B$ is faithfully flat as
 the right $A$-module.  Then we define $M$ as $B\otimes_A B$ with the natural 
structure of a coalgebra in the
 category $B$-$mod$-$B$.  We have a natural morphism of covers $$f:(A,A)\ra 
(B,M)$$

\proclaim {Theorem}.  The morphism $f$ as above defines an equivalence $f^*$ 
between the categories $QCoh(B,M)$
 and $QCoh(A,A)=A$-$mod$.  \par

This theorem follows from the general Barr-Beck theorem. Recall that a {\it 
comonad} (or {\it cotriple}) in the category $C$ consists of a functor 
$T:C\rightarrow C$ and two functor morphisms $\delta :T\rightarrow T^2$ and 
$\epsilon :T\rightarrow Id_C$ such that $T\delta \circ \delta =\delta T\circ 
\delta $ and $\epsilon T\circ \delta =id_T=T\epsilon \circ \delta $. A {\it 
coaction} of the comonad $T$ is a morphism $c:X\rightarrow T(X)$ such that 
$Tc\circ c=\delta (X)\circ c$ and $\epsilon (X)\circ c=id_X$. We denote by 
$T$-$comod$ the category of {\it $T$-comodules}, i.e. objects of $C$ endowed 
with 
a coaction of $T$, with naturally defined morphisms: a morphism 
$(X,c)\rightarrow 
(X',c')$  is a morphism $f:X\rightarrow X'$ such that $Tf\circ c=c'\circ f$. An 
arbitrary pair of adjoint functors $F:C_1\ra C_2$, $G:C_2\ra C_1$ with 
adjunction 
morhisms $\epsilon :FG\rightarrow Id_{C_2}$ and $\eta :Id_{C_1}\rightarrow GF$ 
determines a comonad $(T,\delta ,\epsilon )$, where $T=FG,\ \delta =F\eta G$, 
and 
there is an obvious functor $C_1\ra T$-$comod,\,\,\,, Y\mapsto (F(Y),F\eta 
(Y))$. The 
following version of Barr-Beck's theorem is sufficient for our needs.

\medskip

\noindent {\bf Theorem.} {\it Let $F:C_1\ra C_2$ be a functor between two 
categories having a right adjoint functor
     $G:C_2\ra C_1$.  Assume that the following conditions hold:

(a) a pair of arrows $f,g:X\rightrightarrows Y$ has a kernel if its image by $F$ 
has a kernel;

(b) if $h:Z\ra X$ is such that $f\circ h=g\circ h$ and $Fh$ is the kernel of the 
pair $(Ff,Fg)$, then $h$ is a kernel of $(f,g)$.

Then the canonical functor $C_1\ra T$-$comod$ is an equivalence of  categories.}

\medskip 

We apply this theorem to the case $C_1=A$-$mod$, $C_2=B$-$mod$ and 
$F=B\otimes_A$.

\medskip

Note that the arbitrary pairs of morphisms in Barr-Beck's theorem can be  
replaced 
by so called {\it coreflexive pairs}, i.e. the pairs $f,g:X\rightrightarrows Y$ 
such that there exists a morphism $e:Y\rightarrow X$ such that $e\circ 
f=id_X=e\circ g$. And even this condition can be weakened (cf. [MLM], IV.4).

\medskip

\noindent {\bf 2.1.3.  Refinements of covers}\medskip

Let $(B_1,M_1)$ be a cover and $i:B_1\ra B_2$ be an inclusion of algebras such 
that $B_2$ is faithfully flat as a
 right $B_1$-module.  We define $M_2$ as $B_2\otimes_{B_1} M_1 \otimes_{B_1} 
B_2$ 
with the natural structure of a
 coalgebra in the category $B_2$-$mod$-$B_2$.  It is easy to see that 
$(B_2,M_2)$ 
is again a cover and we have a
 natural morphism of covers $f:(B_2,M_2)\ra (B_1,M_1)$.

\proclaim {Definition}.  A morphism of covers isomorphic to $f$ as above is 
called 
a refinement morphism.  \par

 We denote by $Ref$ the class of refinement morphisms.  Analogously to the 
example 
FC3) from the previous
subsection, the functors $f^*$ for $f\in Ref$ are equivalences of categories of 
quasi-coherent sheaves.  Here
 we apply Barr-Beck theorem to the functor
 which  acts from $QCoh(B_1,M_1)$ to $B_2$-$ mod$ and is given by formula
$F(\E,m_{\E})=B_2\otimes_{B_1}\E$.

The class of refinements is closed under compositions. Also, any pullback of a 
refinement morphism
 is again a refinement morphism. Thus, we can easily describe the localization 
of 
the category
  of covers with respect to refinements. 

\medskip

 \noindent {\bf 2.2.  Noncommutative spaces}\medskip
 
 In this section we give the definition of noncommutative spaces (in several 
steps).
  Morally one should think about noncommutative space $X$ as about  an abelian 
category $QCoh(X)$
 and two adjoint functors $\pi_*:QCoh(X)\ra k$-$mod$, $\pi^*:k$-$mod\ra 
QCoh(X)$. The pullback
 $\pi^*(k^1)$ of the standard 1-dimensional $k$-module is the ``structure 
sheaf'' $\O_X$.
 
 \medskip

 \noindent {\bf 2.2.1.  Covers of noncommutative spaces}\medskip

We define covers of noncommutative spaces adding some data and an 
axiom to the definition of the category of covers:

\proclaim {Definition}.  Objects of category $Covers_k^{sp}$ are pairs $(C,s_C)$ 
where $C=(B,M)$ is a cover and
 $s_C:B\otimes B\ra M$ a homomorphism of bimodules which is an epimorphism and 
also a morphism of coalgebras.
 Morphisms in $Covers_k^{sp}$ are morphisms of covers compatible with structure 
epimorphisms from $B\otimes
 B$.\par

In other words, we demand that $M$ should be a cyclic bimodule generated by an 
element $m:=s_C(1_B\otimes 1_B)$
 which behaves well with
 respect to the coalgebra structure.  Such $M$ is given by a left ideal 
$Ker(s_C)$ 
in $B\otimes
 B^{opp}$ satisfying a complicated system of axioms.

 We will call objects of $Covers_k^{sp}$ {\it space covers}.  Notice that all 
our 
examples of covers (except stacks) are
automatically space covers.

The morphism $s_C$ give rise to a canonical morphism from $(C,s_C)$ to the 
refinement
 $(B,B\otimes B)$ of the  ``point'' object $Spec(k):=((k,k),id)$
(the final object in $Covers_k^{sp}$).  The pullback of the standard 
$1$-dimensional module $k^1$ under this
morphism we call the structure sheaf $\O_C$.  It is represented by 
$1$-dimensional 
free $B$-module $B^1=B$ with
the coaction of $M$ arising from $s_C$.

Also,  for any refinement $f:(B_2,M_2)\ra (B_1,M_1)$ from cover to a space cover
 there exists a unique
structure of space cover on $(B_2,M_2)$ such that $f$ became a morphism of space 
covers.

We  define {\it equivalence classes of noncommutative spaces\rm} as equivalence 
classes of
 space covers modulo the relation generated by refinements. In order define the 
``right''
 {\it category\rm} of noncommutative spaces we need to make an additional work.
 
 \medskip

 \noindent {\bf 2.2.2.  Equivalence relation between morphisms of space 
covers}\medskip
 
 Two different morphisms of space covers could give isomorphic functors
  of pullbacks on categories of quasi-coherent sheaves. We would like to  
identify such morphisms.
  
  Let $(C,s_C)$ and $(C',s_{C'})$ be two space covers, and $B$ resp. $B'$ denote 
  corresponding algebras.
  
  \proclaim {Definition}. Two morphisms $f,g: (C,s_C)\ra (C',s_{C'})$ are 
equivalent iff
   for any element $\sum_{\alpha} x_{\alpha}\otimes y_{\alpha}\in 
Ker(s_{C'})\subset B'\otimes B'$
    the following equations hold:
    $$\sum_{\alpha} f(x_{\alpha})\cdot g(y_{\alpha})=\sum_{\alpha} 
g(x_{\alpha})\cdot f(y_{\alpha})=0\in B$$\par
    
    \proclaim {Theorem}. The relation described  above is an equivalence  
relation. There is a canonical
     identification of pullback functors for equivalent  morphisms.
     \par
     
     We will not give here the proof of this theorem. It follows from general 
categorical considerations
      and from an explicit calculation of certain 
      products in the category of space covers (see theorem in 2.3).
     
 The equivalence of morphisms  is compatible with composition.
 Thus, we can form a quotient category $\widetilde {Covers}_k^{sp}$.
 
\medskip

 \noindent {\bf 2.2.3.  The category of noncommutative spaces}\medskip

\proclaim {Definition}.  The category $Spaces_k$ of noncommutative spaces over 
$k$ 
is defined as the localization
$\widetilde{Covers}_k^{sp}[Ref^{-1}]$ of the category of space covers with 
respect to the 
class of refinements.  \par

We call morphisms of the category $Spaces_k$ ``maps''.  Let us spell explicitly 
how maps look 
like in terms of 
covers. A map from $(B, M)$ to $(B_1, M_1)$ is given by a  refinement $(B',M') 
\ra 
(B_1,M_1)$
 and a morphism of space covers $f':(B',M')\ra (B_1,M_1)$. Two such data  
 $\bigl((B',M'),f'\bigr)$ and $\bigl((B'',M''),f''\bigr)$ give the same map iff 
for some (equivalently, 
  for any) given common refinement $(B''',M''')$ (e.g. the cartesian product of 
$(B',M')$ and $(B'',M'')$
   over $(B,M)$)
   two corresponding  cover morphisms from $(B''',M''')$ to $(B_1,M_1)$ are 
equivalent.

\proclaim {Theorem}.   The category of separated quasi-compact schemes over $k$ 
is 
equivalent to a full
 subcategory of $Spaces_k$.  The category $Alg_k^{opp}$ is also equivalent to a 
full subcategory of $Spaces_k$.
  Category $Spaces_k$ has finite  limits. 
 \par

\proclaim {Theorem}.  Construction $(B,M)\ra QCoh(B,M)$ extends to a
 functor from the category $Spaces_k$ to the category of abelian
 $k$-linear categories. For separated quasi-compact schemes, it gives usual 
quasi-coherent sheaves
  of commutative algebraic geometry. For associative algebras, it gives 
categories 
of left modules. \par

One can
consider noncommutative stacks by eliminating the condition that the structure 
homomorphism
$s_S:B\otimes B\ra M$ is surjective.  We will not develop here the more general 
theory of stacks.

\medskip

\noindent {\bf 2.3.  Affine covers}

\medskip

We call {\it affine spaces\rm} noncommutative spaces corresponding to the cover 
of 
the type $(B,B)$ where $B$ is an
 algebra.  Abusing notations, we will denote the space $(B,B)$ by $Spec(B)$, 
although here
  we do not use points of the spectrum ([R3], [R1], [R4]) at all.  The algebra 
$B$ 
can be identified with 
 $$Hom_{QCoh(Spec(B))}(\O_{Spec(B)},\O_{Spec(B)})\,\,\,.$$
   Also, 
 $B$ coincides as a set with the set of maps from $(B,B)$ to the 
(non)-commutative 
affine line $\A^1_k=N\A_k^1:=Spec(k[t])$.
 In general, $d$-dimensional noncommutative affine space $N\A_k^d$ is defined as 
$Spec(k\langle
 x_1,\dots,x_d\rangle)$.  The algebra structure on $\O(S)=Maps(S,N\A_k^1)$ is 
induced by the algebra structure on $N\A_k^1$ considered as an object of  
$Spaces_k$.

Let $\pi:Spec(B)\ra S$ be a cover of a space by an affine space.  The category 
of 
spaces admits finite limits.
 Thus, we can form an infinite sequence of spaces $S^{(n)},\,\,\,n\ge 1$ taking 
the cartesian product over $S$ of
 $n$ copies of $Spec(B)$.  It is easy to see that all spaces $S^{(n)}$ are 
affine. 
 The algebra of functions on
 $S^{(n)}$ is generated by $n$ copies $i_1(B),\dots, i_n(B)$ of the algebra $B$ 
subject to the set of 
 relations which
 we will describe now.
\medskip

 {\bf Relations:} {\it Let $z=\sum_{\alpha} a_{\alpha}\otimes b_{\alpha}\in 
B\otimes B $ be an arbitrary
 element of the kernel of the structure
map $s_S:B\otimes B\ra M$.  Then for any $k\ne,l,\,\,\, 1\le k,l\le n$ the 
element 
$\sum_{\alpha} i_k(a_{\alpha})
\cdot i_l(a_{\alpha})$ is equal to zero in $\O(S^{(n)})$. \rm}

\medskip

One can also give a description of spaces $S^{(n)}$ in terms of categories of 
quasi-coherent sheaves.

\proclaim{Theorem}. The category $QCoh(S^{(n)})$ is canonically equivalent to 
the
 category of $n$-tuples $(\E_i)_{i=1,\dots, n}$ of objects of 
$QCoh(Spec(B))\simeq B$-$mod$ endowed
 with an identfication of $n$ objects $\pi_*(\E_i),i=1,\dots,n$ in the category 
$QCoh(S)$. \par

It is easy to see that the collection of spaces $S^{(n)}$ (and natural maps 
between them)
 form a contravariant functor
 from the category of finite non-empty sets to the category of affine spaces.  
For 
any  space $U$, the
 set of maps $Maps(U,S^{(2)})$ is a subset of the  square of the set 
$Maps(U,S^{(1)})$ and forms there a graph
 of an equivalence relation.  A map from $U$ to $S$ can be described locally 
(i.e. 
after passing to a cover)
  as a map to
  $S^{(1)}$ modulo the equivalence relation coming from maps to $S^{(2)}$.

\medskip

\noindent {\bf 2.4.  Cohomology of quasi-coherent sheaves}

\medskip

As we mentioned in the previous section, for every space $S$, there exists a 
distinguished object $\O_S$ in the
 category $QCoh(S)$.  It is easy to see that the functor $\E\mapsto 
Hom(\O_S,\E)$ 
is the right adjoint to the
 inverse image functor $f_S^*;QCoh(point)\ra QCoh(S)$ for the canonical map 
$f_S:S\ra point=Spec(k)$.  We denote the vector
 space $Hom(\O_S,\E)$ as $\Gamma(\E)$ (the space of global sections of $\E$).  
The 
functor $\Gamma$ is  left exact.

We would like to define now the derived functor for $\Gamma$.  First of all, it 
makes sense because of the
 following

\proclaim {Theorem}.  Category $QCoh(S)$ has enough injective objects.  

\par

The proof is the following.  Let us chose an affine cover $\pi:  Spec(B)\ra S$.  
The functor $\pi^*$ has a right
adjoint $\pi_*$.  The category $B$-$mod$ (canonically $\simeq QCoh(Spec(B))$) 
has 
enough injective
objects.  Now, for any quasi-coherent sheaf $\E$ on $S$ let us chose an 
embedding 
of $\pi^*\E$ 
into an injective $B$-module
$I$.  Due to the fact that the inverse image functor $\pi ^*$ is exact,  $\pi_* 
I$ 
is again an injective object. Since the functor $\pi ^*$ is faithful, the  
natural 
homomorphisms 
$\E\ra\pi_*\pi^* \E $ is a monomorphism. And the image $\pi_*\pi^* \E \ra
\pi_*I$ of the embedding $\pi^* \E \ra I$ is a monomorphism too, because the 
direct image functor $\pi _*$ is left exact (as any functor having a left 
adjoint). 

\medskip

  Thus, one can proceed and define derived functor $R\Gamma$ (and more generally 
$RHom$, more generally derived functor of any left exact functor) using 
injective
resolutions.  

\proclaim {Theorem}. For any cover $(B,M)$ of $S$  and for any object 
$(\E,m_{\E})\in QCoh(B,M)$, 
 the cohomology of $S$ with
 coefficients in $\E$ can be calculated via the \v Cech 
 complex $$\E\ra M\otimes_{B}\E\ra M\otimes_B M\otimes_B\E\ra
 \dots$$\par 

One can extend the previous result to the relative case.  For every morphism 
$f:(C_1,s_{C_1})\ra (C_2,s_{C_2})$ in
category $Covers_k^{sp}$ the functor $f^*$ admits a right adjoint $f_*$.
 Since {\it the direct image functor} $f_*$ is left exact, there is a derived 
functor for  $f_*$.

\medskip

\noindent {\bf 2.5.  The definition of coherent sheaves}

\medskip

In commutative algebraic geometry the abelian category of coherent modules are 
usually considered in the case of
 noetherian rings.  In the noncommutative setting, our principal examples are 
free 
algebras in several indeterminates which are far from being noetherian.   
Nevertheless, one can define a reasonable abelian category of modules over 
smooth 
algebras as well.

\proclaim {Definition}.  A module $M$ over an algebra $A$ is called coherent iff 
it is finitely presented,
 i.e.  there exists an exact sequence $$F_{-1}\ra F_{0}\ra M\ra 0$$ where 
$F_{-1},\,\,F_0$ are free finitely generated
 $A$-modules.  \par

We denote $Coh(A)$ the full subcategory of $A-mod$ objects of which are coherent 
modules.

\proclaim {Lemma}.  For any hereditary algebra $A$ (e.g. for  a formally smooth 
$A$), 
the category $Coh(A)$ is an abelian category.  \par

It is enough to check that the kernel of any morphism $\phi:M\ra N$ of coherent 
$A$-modules is coherent.
The proof is the following. In hereditary abelian categories any object in the 
bounded derived category
 is isomorphic to its cohomology. Thus, $Ker(\phi)$ is a direct summand in 
$D^b(A-mod)$ of a perfect complex.
 Using telescope construction we see that $Ker(\phi)$ is quasi-isomorphic in an 
infinite complex
  bounded from above of finitely generated projective modules. From this it 
follows that $Ker(\phi)$ is coherent. 

\medskip
   
The property of a module over arbitrary algebra to be (or not to be) finitely 
presented is preserved under
  faithfully flat extensions of the algebra.  This implies that for spaces one 
can 
define the notion of a
  finitely presented quasi-coherent sheaf passing to an arbitrary affine cover.  
The property of being finitely
  presented is independent of the cover.

For any noncommutative space $S$ admitting a cover 
 by a smooth or noetherian affine space, the category of
 finitely-presented quasi-coherent sheaves is abelian.  In such a case we call 
this the category $Coh(S)$ of
 coherent sheaves.
\medskip

\noindent {\bf 2.6.  Formally smooth, formally non-ramified and
 formally \'etale
morphisms.}

\medskip

We call an affine scheme morphism  $Spec(A)\rightarrow Spec(B)$  
a {\it thickenning} if the corresponding algebra morphism $B\rightarrow A$ 
is an epimorphism and its kernel is a nilpotent ideal. 

\medskip

Following the pattern of commutative algebraic geometry [EGA, IV.17], 
we give the following 

\medskip  

\noindent {\bf  Definitions.} {\it 1) Let $f:X\rightarrow Y$ be a  
morphism.  We 
call $f$ formally smooth (resp.   formally non-ramified, resp.
 formally \'etale) if, for any affine scheme $Spec(A)$, 
 any  thickenning $Spec(T)\rightarrow Spec(A)$, and 
any morphism $Spec(A)\rightarrow Y$, the canonical morphism $$
Hom_{Y}(Spec(A),X)\rightarrow Hom_{Y}(Spec(T),X) $$
defined by the immersion $Spec(T)\rightarrow Spec(A)$ is surjective 
(resp.  injective, resp.  bijective). 

2) A $k$-space $X$ will be called {\it formally smooth}   if 
the canonical morphism $X\rightarrow Spec(k)$ is formally smooth. }

\medskip

\noindent {\bf  Remarks.}  (i) If $X$ is affine, then the 
definition of the formal smoothness given here is in
accordance with the one in 1.1.

(ii) To check that $f$ is formally smooth (resp.  formally non-ramified,
resp.  formally \'etale), it suffices to do it in the case when the 
square of the kernel $J$ of the algebra epimorphism $A\rightarrow T$ is zero. 

\medskip

(iii) The properties of $f$ defined in 2.6 are the properties of the
represented by $f$ functor $Y'\mapsto Hom_{Y}(Y',X)$ from the category dual
to the category of $Y$-spaces to the category $Sets$.  They make sense for
any functor $(Spaces_k/Y)^{op}\rightarrow Sets$, representable or not.

(iv)  It follows from the definition that  $f$ is formally \'etale iff it 
is both formally smooth and formally non-ramified. 
One can show that $Spec(A)$ is \'etale 
(i.e. $Spec(A)\rightarrow Spec(k)$ is \'etale)
 iff the algebra $A$ is {\it separable}. 
 The latter means that $A$ is a projective $A$-bimodule, or 
 equivalently, $A$ has dimension zero with respect to Hochschild cohomology. 
 An example of a separable algebra is the algebra $Mat(n\times n,k)$ 
 of $n\times n$ matrices.

\bigskip

\centerline {\bf 3.  \rm NONCOMMUTATIVE PROJECTIVE SPACE}

\medskip

\noindent {\bf 3.1.  Definition \'a la Grothendieck.}

\medskip

For any integer $d\ge 1$, we define the {\it noncommutative projective space 
$N{\P}^{d-1}_k$} as a space
 representing a certain contravariant functor from the category of affine spaces 
to the category $Sets$.

Let $A$ be an algebra.  The set $Map(Spec(A),N{\P}^{d-1}_k)$ should be 
functorially identified with the set of
 quotient modules $\F$ of the standard free $d$-dimensional $A$-module $A^d$, 
such 
that locally in flat topology
 $\F$ is isomorphic to the free $1$-dimensional $A$-module. The last condition 
implies that
 $\F$ is projective.

From this definition, it is not clear whether $N{\P}^{d-1}_k$ exists. We shall 
sketch a prove of its existence by producing an explicit affine cover.

\medskip

\noindent {\bf 3.2.  Explicit cover.}

\medskip

We will construct a cover of $N{\P}^{d-1}_k$ analogous to the Jouanolou's cover 
of 
the projective space in
 commutative algebraic geometry.  Recall that this is a cover  by the affine 
quadric:
 $$\{(x_1,\dots,x_d;y_1,\dots,y_d)|\sum_{i=1}^d x_i y_i =1\}\ra \{(x_1:\dots 
:x_d)\}=\P^d(k)\,\,\,.$$

Let $B$ be the associative algebra generated by $2d$ variables 
$x_1,\dots,x_d,y_1,\dots, y_d$ satisfying the
 relation $\sum_{i=1}^d y_i x_i =1$.  The space $Spec(B)$ represents the 
following 
functor on affine spaces:
 $Maps(Spec(A),Spec(B))$ is the set of quotient modules $\F$ of $A^d$ together 
with an isomorphism $i:A^1\simeq \F$ of
 $A$-modules, and a splitting $s$ of $A^n$ into the direct sum of $\F$ and a 
complementary module $\F'$.  
 Also, for each
 $n\ge 1$ we can consider functor whose value on $Spec(A)$ is the set of  
collections
 $(\F,i_1,s_1,\dots, 
 i_n,s_n))$, 
 where $\F$ is a quotient module of $A^d$ 
and for each $k=1,\dots,n$, the triple
 $(\F,i_k,s_k)$ belongs 
 to $Maps(Spec(A),
 Spec(B))$.  It is easy to see that this functor is representable by an affine 
space whose algebra of functions 
 we
 denote by $B^{(n)}$.

 Thus  we have constructed  a contravariant
  functor from the category of non-empty finite sets to the category of
affine spaces.  We can construct (as in 2.3)
a new contravariant functor from spaces to sets by passing locally
 to sets of equivalence classes of maps to $Spec(B)$ modulo an equivalence 
relation coming from maps
  to $Spec(B^{(2)})$. This new functor
   coincides, more or less by definition, with the functor described in the 
previous section
3.1. 

One can check by direct calculations that the algebra $B^{(2)}$ is generated by 
two copies $i_1(B),\,\,\,i_2(B)$
modulo relations of the type described in 2.3.
  The bimodule $M$ in this case is $(B\otimes B)/I$ where $I$ is the
left $B\otimes B^{opp}$-ideal generated by the following set of elements:  
$$e_j:=- 1\otimes x_j
+\sum_{i=1}^d x_j y_i \otimes x_i
,\,\,\,\,j=1,\dots,d\,\,.$$

This shows that the functor defined in 3.1 is representable by $(B,M)$.

\medskip

Analogously, one can imitate the Grothendieck's definition of projective space 
for 
coherent sheaves and introduce
relative projective space $\P(\E)$ for any noncommutative space $S$ and a 
finitely 
generated quasi-coherent sheaf
$\E$ on $S$.

\medskip

\noindent {\bf 3.3.  Derived category of quasi-coherent sheaves}

\medskip

Here we will only state the results without proofs (proofs will appear among 
other 
things in [KR]).  We denote by $Q_d$ the quiver which has two vertices
$\{v_0,v_1\}$ and $d$ oriented edges all going from $v_0$ to $v_1$.

\proclaim {Theorem}.  For any $d\ge 1$ the category $QCoh(N\P^{d-1}_k)$ has 
cohomological dimension $1$.  The bounded
derived category $D^b(QCoh((N\P^{d-1}_k))$ is equivalent to the bounded derived 
category of representations of the
quiver $Q_d$.  \par

\proclaim {Theorem}. The algebra $B$ described in 3.2  and the projective space 
$N\P^{d-1}_k$ are formally smooth.  The 
bounded derived category of
coherent sheaves on $N\P^{d-1}_k$ is of finite type (i.e.  for any two objects
$\E,\F$ we have $\sum_i rk(RHom^i(\E,\F))<+\infty$), and it is equivalent to the 
bounded derived category of 
finite-dimensional representations of $Q_d$.
The  group $K_0(Coh(N\P^{d-1}_k))$ is  free abelian group with $2$ generators.  
\par

The main ingredient in the proof of both theorems is a noncommutative analogue 
of 
Beilinson's resolution of the
sheaf of functions on the diagonal in $\P^{d-1}_k\times \P^{d-1}_k$.  The 
noncommutative resolution is much shorter
 than the commutative one, and has length $2$.  There are two coherent sheaves 
$\O$ and  $\O(1)$ 
 (the universal quotient module $\F$ from 3.1) which generate in an appropriate
sense the whole derived category on $N\P^{d-1}_k$.

There are two surprising particular cases:  $d=1$ for which the noncommutative 
projective space is more
 complicated than the usual one (which is a point), 
and the case $d=2$.  It is well-known that
  the category $D^b(Coh(\P^1_k))$ is
 equivalent to $D^b(Q_2-mod_{finite})$.  Thus, in the case $d=2$ we have three 
different
 abelian categories with the same derived category:  the noncommutative 
projective 
line, the 
 commutative projective line,
 and the quiver $Q_2$.

\medskip

\noindent {\bf 3.4.  Grassmanians}

\medskip

One can define Grassmanians $NGr_k(d',d)$ for $d',d\ge 1$ in the same fashion as 
the projective space, changing the
 requirement that the quotient module $\E$ of $\O^d$ should be locally 
isomorphic 
to a free $d'$-dimensional module,
 instead of a free $1$-dimensional module.  We claim that in contrast to the 
commutative case we do not get a new
 space.  It follows from the fact that there exists a {\it non-zero\rm} algebra 
$A$ such that $A^1$ is isomorphic 
 $A^2$ (and
 thus to $A^{d'}$) as $A$-module.  For example, one can take 
  $A$ equal to $End(V)$ where
  $V$ is an infinite-dimensional vector space.  On the faithfully flat
 extension which is obtained by taking free product with such algebra $A$, we 
identify conditions on $\E$ for all
 Grassmanians.  Thus, $NGr_k(d',d)$ coincides with $N\P^{d-1}_k=NGr_k(1,d)$.

  We expect that various remarkable identities in
   noncommutative linear algebra in free skew field discovered by Gelfand and 
Retakh (see [GR])
    can be interpreted  as identities between morphisms of coherent sheaves on 
$N\P^{d-1}_k$
     and on similar spaces.

\medskip

\noindent {\bf 3.5.  Representation spaces}

\medskip

For a noncommutative space $S$ we can try to define representation schemes 
$Repr_n^S$. At least, it is clear how  sets of $k$-points  should look like:
 $$Reps_n^S(k)=Maps(Spec(Mat(n\times n,k)),S)\,\,\,.$$
 It follows from the definition of the projective space in 3.1 and from the 
discussion in 3.4 that
 the set of $k$-points of $N\P^{d-1}_k$ is equal to
 the set of {\it all\rm} non-zero quotient spaces of $k^d$. We leave to the 
reader 
as an exercise
  the description of maps from  $ Spec(Mat(n\times n,k))$ to $N\P^{d-1}_k$
  for the case $n\ge 2$.

  Let us express again our hopes: representation schemes $Repr_n^S$ for 
(formally)
   smooth 
non-affine spaces $S$ should carry all structures
   described in Part 1 (the double tangent bundle, Kapranov's formal 
noncommutative functions etc.).

\medskip 

\bigskip

\bigskip

\centerline{{\bf References.}  }

\medskip

[CQ1] J. Cuntz, D. Quillen, Algebra extensions and nonsingularity, Journal of 
AMS, 
v.8, no. 2, 1995, 251--289

[CQ2] J. Cuntz, D. Quillen, Cyclic homology and nonsingularity, Journal of AMS, 
v.8, no. 2, 1995, 373--442

[EGA] A.  Grothendieck, J.A.  Dieudonn\'e, El\'ements de G\'eom\'etrie
Alg\'ebrique, Springer Verlag, New York - Heidelberg - Berlin, 1971

[GZ] P.  Gabriel and M.  Zisman, Calculus of fractions and homotopy theory,
Springer Verlag, Berlin-Heidelberg-New York, 1967

[GR] I.  Gelfand and V.   Retakh,  Quasideterminants, I,  Selecta Math., v.3, 
1997, 517--546

[K] M. Kontsevich, Formal non-commutative symplectic geometry, in "The Gelfand 
Mathematical Seminars, 1990-1992", Birkh\"auser, Boston-Basel-Berlin (1993), pp. 
173--187

[Ka] M.  Kapranov, Noncommutative geometry based on commutator expansions,

\noindent math.AG/9802041 (1998), 48 pp.

[Kn] D. Knutson, Algebraic spaces, LNM 203, Springer-Verlag, 1971

[KR] M. Kontsevich, A. Rosenberg,  Noncommutative separated spaces, in 
preparation

[ML] S.  Mac-Lane, Categories for the working mathematicians, Springer -
Verlag; New York - Heidelberg - Berlin (1971)

[MLM] S. Mac Lane, I. Moerdijk, Sheaves in Geometry and Logic, Springer-Verlag,  
New York - Heidelberg - Berlin-London-Paris (1992)

[R1] A.L.  Rosenberg, Noncommutative algebraic geometry and representations of
quantized algebras, Kluwer Academic Publishers, Mathematics and Its
Applications, v.330 (1995), 328 pages.

[R2] A.L.  Rosenberg, Noncommutative schemes, Compositio Mathematica 112 (1998), 
93-125

[R3] A.L.  Rosenberg, Noncommutative local algebra, Geometric and Functional
Analysis (GAFA), v.4, no.5 (1994), 545-585

[R4] A.L.  Rosenberg, The spectrum of abelian categories and reconstruction of
schemes, in "Algebraic and Geometric Methods in Ring Theory", Marcel Dekker,
Inc., New York, (1998), 255-274

\end